\documentclass[12pt]{amsart}
\usepackage[active]{srcltx}
\usepackage{calc,amssymb,amsthm,amsmath,amscd, eucal,ulem}
\usepackage{alltt}
\RequirePackage[dvipsnames,usenames]{color}

\input{mabliautoref.sty}
\input{xy}
\xyoption{all}
\input{kmacros3.sty}
\usepackage{tikz}
\usepackage{amsfonts, mathrsfs}
\usepackage{calligra}
\usepackage{stmaryrd} 
\usepackage[left=1.1in,top=1in,right=1.1in,bottom=1in]{geometry}

\theoremstyle{plain}
\newtheorem*{Question}{Question}

\DeclareMathOperator{\mSpec}{m-Spec}

\newcommand{\bit}{\begin{itemize}}
\newcommand{\eit}{\end{itemize}}
\newcommand{\ben}{\begin{enumerate}}
\newcommand{\een}{\end{enumerate}}
\newcommand{\bpf}{\begin{proof}}
\newcommand{\epf}{\end{proof}}

\newcommand{\OO}{\mathscr{O}}

\title{Counterexamples to Bertini Theorems for Test Ideals}
\author{Andrew Bydlon}
\thanks{The author is grateful to have been supported by the NSF FRG Grant DMS-\#1501115, as well as the NSF RTG Grant at the University of Utah, DMS \#1246989}
\begin{document}
\begin{abstract}
In algebraic geometry, Bertini theorems are an extremely important tool.  A generalization of the classical theorem to multiplier ideals show that multiplier ideals restrict to a general hyperplane section.  The test ideal can be seen to be the characteristic $p > 0$ analog of the multiplier ideal.  However, in this paper it is shown that the same type of Bertini type theorem does not hold for test ideals.
\end{abstract}
\maketitle
\section{Introduction}\label{Sect:Intro}

In characteristic $0$ algebraic geometry, Bertini theorems are often used to prove theorems using the idea of descending induction.  That is to say if a variety exhibits a given property $\mathcal{P}$, it is often desirable for the property to descend to a general choice of codimension 1 subvariety.  Many properties in characteristic $0$, such as being smooth or having relatively nice singularities, do descend to a general hyperplane section.  More generally, it was later shown that an important invariant, the multiplier ideal, also exhibits these Bertini type theorems for pairs (or even triples) $(X,\Delta)$.

The case of characteristic $p>0$ was less understood.  First, the notion of a multiplier ideal is not known to exist unless a log resolution of the pair $(X,\Delta)$ exists.  Furthermore, the positive characteristic analog of the multiplier ideal, the test ideal, was not known to satisfy such a property except in very specific circumstances such as the case where $\tau(X,\Delta) = \mathscr{O}_X$.  In this paper, it is shown in \autoref{prop:FailSB1} that the natural analog of the Bertini theorem for multiplier ideals is false for test ideals in dimension $n\geq 3$.

\begin{theorem*}
There exist $\mathbb{Q}$-divisors $\Delta$ on $\PP^n$ with the property that for a general hyperplane $H$,
$$\tau(\PP^n,\Delta)|_H \neq \tau(H,\Delta|_H).$$
\end{theorem*}

In the process of proving this proposition, it is shown in \autoref{thm:Cex} that the analog of Koll\'ar-Bertini theorem \cite[9.2.29]{Laz} is also false for a class of pairs.

\begin{theorem*}
	Given any $\epsilon>0$, and $n\geq 3$, there are pairs $(\PP^n,\Delta)$ such that for a general hyperplane $H$,
	$$\tau(\PP^n,\Delta) \neq \tau(\PP^n,\Delta+\epsilon H).$$
\end{theorem*}

This in some sense eliminates any hope that a Bertini theorem can hold for test ideals. Thus for many types of F-singularities which are governed by the test ideal, taking a general hyperplane section will not preserve the type of singularity.  However, it is shown to hold in \cite{SZ} for pairs $(X,\Delta)$ which are sharply F-pure or strongly F-regular.  In \autoref{Sect:Questions}, some follow-up questions are posed regarding the open state of the question for non-pairs, the property of being F-rational, and whether there is something more geometric to consider regarding these types of questions. 

\

\noindent{\it Acknowledgements}: I would like to thank my advisor Karl Schwede for suggesting this problem, as well as many valuable discussions that led to understanding the nature of the problem as well its the eventual solution.  His assistance with earlier drafts of this paper is also greatly appreciated.  The author is also very grateful for the comments of the referee which improved the paper through very detailed feedback.  I also would like to thank Stefano Filipazzi for discussions about many of the corresponding statements for multiplier ideals.  Finally, I am very appreciative for some helpful comments from Eric Canton on a earlier draft of this paper.

\section{Background and Motivation}\label{Sect:Motiv}

\subsection{Test Ideals}

Let $X=\Spec(R)$ be a normal variety over a perfect field $K$ of characteristic $p>0$ and let $\Delta$ be an effective $\mathbb{Q}$-divisor on $X$.  This is the data of a pair $(X,\Delta)$.  The (big) test ideal associated to $R$, denoted $\tau(X,\Delta)$, is defined as follows:

\begin{definition}
	$\tau(X,\Delta)$ is the smallest non-zero ideal $J$ of $R$ such that for every $e>0$ and every $R$-linear homomorphism $\phi: F_*^eR(\lceil (p^e-1)\Delta \rceil) \rightarrow R$, $J$ is compatible with $\phi$ in the following sense:
	$$\phi(F_*^e J)\subseteq J$$
\end{definition}

This ideal can be computed locally and can thus be extended to non-affine varieties.  It is a very important invariant which governs the singularities of a positive characteristic variety.  In general, $X$ is strongly F-regular (in particular normal and Cohen-Macaulay) if and only if $\tau(X) = \mathscr{O}_X$, and smaller test ideals imply more severe singularities.  In the case that $\Delta$ is a $\Q$-Cartier divisor, and $\Delta=t\cdot\Div(f)$, for simplicity one often writes $\tau(R,\Delta) = \tau(R,f^{t})$. 

\subsection{Relating Multiplier Ideals to Test Ideals}

Let $X \subseteq \PP^n$ be a smooth quasi-projective variety over $\C$.

\begin{definition}\label{DefnMult}
Let $\Delta\geq 0$ be a $\mathbb{Q}$-divisor such that $K_X+\Delta$ is $\Q$-Cartier, and let $\pi:X' \rightarrow X$ be a log resolution of $\Delta$.  The multiplier ideal sheaf is defined as 
$$\mathcal{J}(X,\Delta) = \pi_* \mathscr{O}_{X'}(\lceil K_{X'} - \pi^*(K_X + \Delta)\rceil).$$
\end{definition}

This notion is independent of the log resolution chosen, and measures the singularities of the pair $(X,\Delta)$ in characteristic zero.  The primary motivating theorem is as follows (see \cite[9.2.29 and 9.5.9]{Laz}):

\begin{theorem}\label{KolBer}
Let $(X,\Delta)$ be a pair, with $X$ a smooth complex variety and $\Delta\geq 0$ a $\mathbb{Q}$-divisor.  Let $|V|$ be a base-point free linear series on $X$.  Then for a general choice of divisor $B\in |V|$, and any $0 \leq \epsilon < 1$,
$$\mathcal{J}(X,\Delta) = \mathcal{J}(X,\Delta+\epsilon B).$$ 

In addition, for general $B\in |V|$
$$\mathcal{J}(X,\Delta) \otimes_{\mathscr{O}_X}\mathscr{O}_B \cong \mathcal{J}(B,\Delta|_B).$$
\end{theorem}

So multiplier ideals behave as well as one could expect when taking general hyperplane sections.  Moreover, multiplier ideals are intimately related to test ideals.  Let $X = \Spec(R)$ be a complex variety, with $R=\C[x_0,\ldots,x_n]/I$ with $D\geq 0$ an effective divisor.  If $I=\langle f_1,\ldots,f_m\rangle$, then let $A$ be the $\Z$-algebra generated by all of the coefficients of the $f_i$ and the equations of a log resolution $\pi: \tilde{X} \rightarrow X$.  Then $X$ and $\pi$ are defined over $A$.  For each maximal ideal $\mathfrak{m}\subseteq A$, $A/\mathfrak{m}$ is a finite field, and we can consider $X_\mathfrak{m}$ obtained by base change.  The following theorem motivates the idea that many of the well studied properties of multiplier ideals have a chance to descend to test ideals:

\begin{theorem}\label{thm:Red}
\cite{Tak}, cf. \cite{Smi},\cite{HY},\cite{Hara} Let $X$ be a complex projective variety, and $\Delta$ a $\mathbb{Q}$-divisor on $X$ with $K_X+\Delta$ a $\Q$-Cartier divisor.  Then for a given model $A$ of $X$ and an open dense collection of closed points $\mathfrak{m} \in \mSpec(A)$

$$\mathcal{J}(X,\Delta)_\mathfrak{m} = \tau(X_\mathfrak{m},\Delta_\mathfrak{m}).$$
\end{theorem}

In particular, it seems reasonable to ask if the same Bertini-type theorem holds for test ideals as well.  Note that in characteristic $p>0$, one typically needs to replace the consideration of a base point free linear series with the assumption that the condition holds on fibers of a map.  A particular case of this more general conjecture is known in the case of a strongly F-regular pair.  This corresponds to the case of $\tau(X,\Delta) = \mathscr{O}_X$, and the theorem is restated here:

\begin{theorem}\label{thm:SZ}
	\cite[Theorem 6.7]{SZ} Suppose that $X$ is a variety over an algebraically closed field $k$, let $\Delta \geq 0$ be a $\mathbb{Q}$-divisor on $X$, and let $\phi:X \rightarrow \PP^n_k$ be a $k$-morphism with separably generated residue field extensions\footnote{Note that one cannot expect the same statement to hold for a general member of an arbitrary basepoint free linear system.}.  Suppose $(X,\Delta)$ is a strongly F-regular pair.  Then for an open dense subset $U\subseteq (\PP^n_k)^\vee$ and all $H\in U$, $\left(\phi^{-1}(H),\Delta|_{\phi^{-1}(H)}\right)$ is strongly F-regular.
\end{theorem}

Furthermore, one can easily test many lower dimensional hypersurface examples of a fixed degree in a finite field, using Macaulay2 \cite{M2}.  There is a script created by the author available as \cite{M2Code}, which produces a list of hyperplanes for which the restriction fails.  Tests examples included hypersurface degree smaller than $6$ in characteristics $p\leq 11$ (as well as small field extensions) in $\P^n$ for $n\leq 4$.  This yielded only examples for which a few special problematic hyperplanes for restriction theorems (can be excluded as a component of a closed set upon taking the perfection of the underlying field).

\section{Computation of Test Ideals for $(\mathbb{A}^n,f^{\frac{1}{p}})$}\label{Sect:Comp}

Consider the pair $(\mathbb{A}_K^n,f^{\frac{1}{p}})$.  As noted above $f^{\frac{1}{p}}$ represents the divisor defined by $f$ with coefficient $\frac{1}{p}$, and $f$ is an element of the ring $S=K[x_1,\ldots,x_n]$.  Assume $K$ is an $F$-finite field with a basis for $F_*K$ over $K$ given as $\langle F_*u_1,\ldots,F_*u_m \rangle_K$.  Let
$$F_* f=\sum_{\boldsymbol{\alpha},i} s_{\boldsymbol{\alpha},i} F_* u_i \boldsymbol{x^\alpha}.$$
In general, $\boldsymbol{\alpha}$ will denote a multi-index in $\Lambda = \Lambda_{p,n} := \left\{0,1,\ldots,p-1\right\}^{\oplus n}$, coming from the set of exponents occurring in the standard basis of monomials $F_*S$ over $S$.  Similarly, $\boldsymbol{x^\alpha}$ is shorthand for $x_1^{\alpha_1}\cdots x_n^{\alpha_n}$.  The index $i$ will range from $1$ to $m$, to indicate the elements of the basis of $F_*K$ as above.  In general, a boldface symbol will denote a multi-index. Also note that $s_{\boldsymbol{\alpha},i}\in S$ is viewed naturally as a subset of $F_*S$ under the inclusion $S\hookrightarrow F_*S: x \mapsto xF_*1 = F_*x^p$.  Let $\Phi$ denote the generator of $\Hom_S(F_*S,S)$ as an $F_*S$-module under pre-multiplication
$$\Phi: F_*S \rightarrow S: F_* u_1\boldsymbol{x^{p-1}} \mapsto 1$$
\noindent and $\Phi$ sends all of the other elements of the standard basis $\langle F_*u_i\boldsymbol{x^\alpha} \ | \ \boldsymbol{\alpha} \in \Lambda,i=1,\ldots,m\rangle_S$ to $0$.

  The test ideal $\tau(\mathbb{A}^n,f^{\frac{1}{p}})$ can then be realized as the image of $\Phi(F_*f\cdot -)$ which can be computed explicitly as follows.  This was written down as early as \cite{HT}, but the proof following the notation outlined above is provided here for convenience:
  
\begin{proposition}\label{prop:Comp}
	\cite{BMS} $\tau(\mathbb{A}^n,f^{\frac{1}{p}}) = \Phi(\langle F_*  f \rangle_{F_*S}) = \langle s_{\alpha,i} \ | \ \boldsymbol{\alpha} \in \Lambda,i=1,\ldots,m\rangle_S$.
\end{proposition}  
  \bpf
  The first equality follows by definition and the fact that $\Phi^{\circ e}$ generates $\Hom_S(F_*^e S,S)$ for each $e>0$, so it only remains to show that $\Phi(\langle F_*  f \rangle) = \langle s_{\alpha,i}  \ | \ \boldsymbol{\alpha},i\rangle_S$.  Consider an element $g$ of $\Phi(\langle F_*  f \rangle_{F_*S})$.  Then $S$-linearity of $\Phi$ implies
	$$g=\Phi(F_*h\cdot f) = \sum_{\boldsymbol{\alpha},m} s_{\boldsymbol{\alpha},i} \Phi(F_* (u_i \cdot h\cdot \boldsymbol{x^\alpha}))\subseteq \langle s_{\alpha,i}  \ | \ \boldsymbol{\alpha},i\rangle_S.$$
	
	\noindent Similarly, considering the element $\Phi\left(F_* (u_i^{-1}\cdot\boldsymbol{x^{p-1-\alpha}} \cdot f)\right) = s_{\alpha,i} \in \Phi(\langle F_* f \rangle_{F_*S})$ yields the reverse inclusion.\epf

Let $H$ be a hyperplane section in $\mathbb{A}^n$ defined by $l=c_0 + c_1x_1 + \ldots + c_nx_n$. Considering the pair $(\mathbb{A}^n_K,(f\cdot l)^{\frac{1}{p}})$ where we can conclude identically to the previous case that the corresponding test ideal $\tau(\mathbb{A}^n_K,(f\cdot l)^{\frac{1}{p}})$ is given as the image of $\Phi(F_* (f\cdot l \cdot S))$, or equivalently
$$\left\langle c_{i,0} s_{\boldsymbol{\alpha},i} + \sum_{j=1,\ldots,n} c_{i,j} \cdot s_{\boldsymbol{\alpha}-\boldsymbol{1}_j,i} \cdot x_j^{\left\lfloor\left(\frac{\boldsymbol{(\boldsymbol{\alpha}-\boldsymbol{1}_j)+\boldsymbol{1}_j}}{p}\right)_j\right\rfloor} \ | \ \boldsymbol{\alpha}\in\Lambda,i=1,\ldots,m \right\rangle_S$$ 
Note that $\boldsymbol{1}_j$ is representative of a multi-index with $1$ in position $j$ and $0$ elsewhere. Moreover when $\alpha_j=0$, $\boldsymbol{\alpha}-\boldsymbol{1}_j$ is treated as $p-1$ in the $j^{\text{th}}$ entry.  Therefore, if $\alpha_j=0$, then the exponent is treated as
$$\left\lfloor\left(\frac{(\boldsymbol{\alpha}-\boldsymbol{1}_j)+\boldsymbol{1}_j}{p}\right)_j\right\rfloor = \left\lfloor\frac{(p-1) + 1}{p}\right\rfloor =1.$$
Finally, $c_{i,j}\in K$ are elements s.t. $F_*(c_j) = \sum_{i} c_{i,j}F_*(u_i)$.

So rephrasing the statement from above, to disprove the direct analog of \autoref{KolBer}, it is enough to show that for a non-closed (or even open) set of $H=V(l)$, the two ideals $\Phi(F_*f\cdot S)$ and $\Phi(F_*f\cdot l\cdot S)$ do not agree.  This will be shown below as \autoref{thm:Cex}.

One can further consider the case of the test ideal sheaf $\tau(\PP^n,f^{\frac{1}{p}})$, which on each element of the standard affine cover is $\tau(\mathbb{A}^n,f^{\frac{1}{p}}|_{x_i=1})$ since the big test ideal is defined locally.  This is recorded in \autoref{cor:proj} and \autoref{thm:HHfalse} for the sake of completeness.

\section{Examples and Counterexamples}\label{Sect:CountEx}

Let $K=K^{p}$ be a perfect field, and $S = K[x_1,\ldots,x_n]$, and consider $\mathbb{A}^n = \Spec(S)$ in which case the trace map can be view as
$$\Phi: F_*S \rightarrow S.$$

Continuing the use of the multi-index notation $\boldsymbol{x^\alpha} = x_1^{\alpha_1}\cdots x_n^{\alpha_n}$, $\Phi$ is defined by

\[ \Phi(\boldsymbol{x^\alpha}) = \begin{cases} 
      \boldsymbol{x^{\frac{\alpha-p+1}{p}}} & p \ | \ \alpha_i+1 \ \ \forall i \\
      0 & \text{otherwise}
   \end{cases}
\]

The best case scenario for extending Bertini theorems to characteristic $p>0$ would be a direct generalization of the Bertini Theorem for Multiplier Ideals.  That is to say that for $I$ homogeneous, $R = S/I$, and $\Delta$ a $\mathbb{Q}$-Cartier Divisor on $X=\Spec(R)$, is it true that for $0\leq\epsilon<1$ and a general hyperplane $H$, the analog of \cite[9.2.29]{Laz} holds true?
$$\tau(R,\Delta) \stackrel{?}{=} \tau(R,\Delta+\epsilon H)$$
$$\tau(X,\Delta)|_{H} \otimes_{\mathscr{O}_X} \mathscr{O}_H \stackrel{?}{=} \tau(H,\Delta|_H)$$

The following examples shows that this is not the case.  All explicit counter-examples are worked out in the affine case $S=K[x_1,\ldots, x_n]$.  

\begin{theorem}\label{thm:Cex}
	Suppose that $S = K[x_1,\ldots,x_n]$ is a polynomial ring over an infinite perfect field of characteristic $p>0$, with $n\geq 3$.  Consider an element $f\in S$ of the form
	$$F_*f = f_{p-1}\cdot F_*(x_1\cdots x_{n-1})^{p-1}x_n^{p-1} + f_{p-2}\cdot F_*(x_1\cdots x_{n-1})^{p-1}x_n^{p-2} + \ldots + f_0\cdot F_*(x_1\cdots x_{n-1})^{p-1}$$
	
	\noindent where the $f_i$ are polynomials in $K[x_1,\ldots,x_{n-1}]$ with the additional independence property that for each $i$,
	\begin{equation}f_i \notin \langle f_0,\ldots,f_{i-1}\rangle_S+ \mathfrak{p}\cdot \tau(S,f^{\frac{1}{p}}) = \langle f_0,\ldots,f_{i-1},x_j f_i,\ldots, x_jf_{p-1}\rangle_S \tag{$\star$}\end{equation}
	
	\noindent where $\mathfrak{p}=\langle x_1,\ldots,x_{n-1}\rangle$ and $j=1,\ldots,n-1$.  Then for a general hyperplane $H=V(l)$,
	$$\tau(S,f^{\frac{1}{p}}) \neq \tau(S,(l\cdot f)^{\frac{1}{p}}).$$
\end{theorem}
\bpf
	Following \autoref{prop:Comp}, one can conclude directly that $\tau(S,f^{\frac{1}{p}}) = \langle f_0,\ldots,f_{p-1}\rangle_S$ and giving the linear form $l$ a presentation of $l=c_0+c_1x_1 + \ldots + c_nx_n$, we can conclude $\tau(S,(f\cdot l)^{\frac{1}{p}})$ is generated by elements of the form:
	
\	
	
	\begin{itemize}
		\item $c_0f_{i} + c_nf_{i-1}$ for $i=1,\ldots,p-1$, the coefficient of $F_*(x_1\cdots x_{n-1})^{p-1}x_n^i$.
		\item $c_0f_0 + c_nx_nf_{p-1}$, the coefficient of $F_*(x_1\cdots x_{n-1})^{p-1}$.
		\item $c_jx_jf_i$ for $i=0,\ldots,p-1$ and $j=1,\ldots,n-1$, the coefficient of $F_*(x_1\cdots \hat{x}_j\cdots x_{n-1})^{p-1}x_n^i$.
	\end{itemize}

\

	Restrict to considering only those $l$ for which all of the $c_j$ are non-zero, which will stand as the open condition in the theorem.  To prove the claim, note that it is enough to show that $f_{p-1} \notin \tau(S,(f\cdot l)^{\frac{1}{p}})$, which by considering the generators of the first type above is equivalent to showing every $f_i\notin \tau(S,(f\cdot l)^{\frac{1}{p}})$. Suppose (aiming for a contradiction) $f_{p-1}$ has a presentation as an element of $\tau(S,(f\cdot l)^{\frac{1}{p}})$:
	$$f_{p-1} = g_0\cdot(c_0f_0 + c_nx_nf_{p-1}) + \sum_{i\leq p-1} g_i\cdot(c_0f_{i} + c_nf_{i-1}) + \sum_{i,j} h_{i,j} \cdot x_jf_i.$$
	
	\noindent where $g_i,h_{i,j}\in S$ are some coefficients, and the indices $i,j$ come from the list of generators above. Without loss of generality, we may assume that $g_i\in K[x_n]$ by taking any monomial terms of the original presentation of $g_i$ in $\langle x_1,\ldots,x_{n-1}\rangle_S$, and instead include them in some combination of the $h_{i,j}$. Explicitly, if each $g_i = g_{i,0} + x_1g_{i,1} + \ldots + x_{n-1}g_{i,n-1}$ with $g_{i,0}\in K[x_n]$, then one can replace $g_i$ with $g_{i,0}$, replace $h_{p-1,j}$ with $h_{p-1,j} + c_0g_{p-1,j}+c_nx_ng_{0,j}$, and replace each $h_{i,j}$ with $h_{i,j}+c_0g_{i,j}+c_ng_{i+1,j}$ for each $i<p-1$.  Rearranging the original representation of $f_{p-1}$ yields  
	\begin{equation}(1-c_nx_ng_0-c_0g_{p-1})f_{p-1} = \sum_{i\leq p-2} (c_0g_i+c_ng_{i+1})f_i + \sum_{i,j} h_{i,j} \cdot x_jf_i.\tag{$\dagger$}\end{equation}
	
	Since each $f_i$ are polynomials only in the variable $x_1,\ldots,x_{n-1}$, the condition $(\star)$ implies that $x_n^af_i \notin x_n^a\cdot\langle f_0,\ldots,f_{i-1}\rangle_S+ x_n^a\cdot\mathfrak{p}\cdot \tau(S,f^{\frac{1}{p}})$ for any $a>0$.  This allows one to consider $(\dagger)$ and the equations that follow filtered by their $x_n$-degree.  Indeed, suppose \linebreak $1-c_nx_ng_0-c_0g_{p-1} \neq 0$ and consider the smallest $x_n$ degree of $1-c_nx_ng_0-c_0g_{p-1}$, and write it as $Cx_n^a$ with $C\in K^\times$.  Then considering the whole of equation $(\dagger)$ at $x_n$ degree $a$,
	$$Cf_{p-1} = C_0f_0 + \ldots + C_{p-2}f_{p-2} + \sum_{i,j} h_{i,j}'x_jf_i.$$
	
	\noindent with $C_i\in K$, and $h_{i,j}'\in K[x_1,\ldots,x_{n-1}]$ is the $x_n$-degree $a$ component of $h_{i,j}$.  However, ($\star$) implies directly that this is not possible, providing a contradiction to the fact that $C\neq 0$, which is to say that there is no lowest $x_n$-degree component of $1-c_nx_ng_0-c_0g_{p-1}$.  Thus, $g_{p-1} = c_0^{-1}(1-c_nx_n g_0)$.
	
	\begin{claim}\label{thm:claim}
	The condition $(\star)$ implies that $c_ng_{i} = -c_0g_{i-1}$ for each $i=1,\ldots,p-1$.
	\end{claim}
	
	\bpf
	The claim is proved by descending induction on $i$, starting with the case $i=p-1$.  Rearranging the presentation above again using the fact that $1-c_nx_ng_0-c_0g_{p-1}=0$ yields	
	$$-(c_ng_{p-1}+c_0g_{p-2})f_{p-2} = \sum_{i\leq p-3} (c_0g_i+c_ng_{i+1})f_i + \sum_{i,j} h_{i,j} \cdot x_jf_i$$
	
\noindent The same technique used on ($\dagger$) above, by considering the lowest $x_n$-degree piece, implies that $c_ng_{p-1}+c_0g_{p-2}=0$.	Now, assume $k\geq 1$ and the claim is true for $k+1,\ldots,p-1$.  Then
	
	$$-(c_ng_{k+1}+c_0g_{k})f_{k} = \sum_{i\leq k-1} (c_0g_i+c_ng_{i+1})f_i + \sum_{i,j} h_{i,j} \cdot x_jf_i.$$
	
\noindent The same argument used in ($\dagger$) again shows that $c_ng_{k+1}+c_0g_k = 0$, or equivalently $g_{k+1} = -c_n^{-1}c_0g_{k}$.\epf
	
	Combining all of the data provided by \autoref{thm:claim}, one sees that
	\begin{equation}c_0^{-1}(1-c_nx_n g_0) = g_{p-1} = (-c_0c_n^{-1})g_{p-2} = (-c_0c_n^{-1})^2g_{p-3} = \ldots = (-c_0c_n^{-1})^{p-1}g_{0}.\tag{$\dagger\dagger$}\end{equation}
	
\noindent On one hand, this implies in particular that no $g_i$ is $0$ as that would imply that all $g_i$ were $0$ and the left most side of ($\dagger\dagger$) would read $c_0^{-1} = 0$.  On the other hand, ($\dagger\dagger$) is impossible since the $x_n$-degree of the left hand side of ($\dagger\dagger$) is exactly one larger than that of the right hand side.  Thus no such presentation can exist, and 
	$$f_{p-1} \in \tau(S,f^{\frac{1}{p}}) \setminus \tau(S,(l\cdot f)^{\frac{1}{p}}).$$\epf
	
\begin{remark}
This proof can easily be extended to the case where $K$ is any F-finite infinite field and $F_*^ef\in F_*^eS$, and the conclusion being that for a general choice of hyperplane $H=V(l)$, $\tau(\mathbb{A}^n,f^{\frac{1}{p^e}})\neq  \tau(\mathbb{A}^n,(l\cdot f)^{\frac{1}{p^e}})$, with nearly identical  assumptions and proof.  Namely, choose 
$$f = \sum_{i,j} f_{i,j}\cdot F_*^e(u_j x_1\cdots x_{n-1})^{p^e-1}x_n^i$$

\noindent meeting the assumption
$$f_{i,j} \notin \langle f_{0,k},\ldots,f_{i-1,k},f_{i,1},\ldots,\hat{f}_{i,j},\ldots,f_{i,m}\rangle_S+ \mathfrak{p}\cdot \tau(S,f^{\frac{1}{p^e}})$$

\noindent here, $\mathfrak{p} = \langle x_1,\ldots,x_{n-1}\rangle_S$, $i=0,\ldots,p^e-1$, and $j,k=1,\ldots,m$ represent the basis of $F_*K$ over $K$.  Thus in dimension greater than $2$, this provides examples of pairs with the property that for a given $t>0$, $\tau(S,f^{\frac{1}{p^e}}) \neq \tau(S,l^t\cdot f^{\frac{1}{p^e}})$, simply by taking $e\gg 0$ such that $\frac{1}{p^e} \leq t$.
\end{remark}

This gives a somewhat algorithmic way to produce counterexamples in dimension greater than 2.  Indeed, once the dimension is 3 or larger one can always find $f_i$ meeting the condition $(\star)$ in the theorem.

Additionally, one can homogenize the equation for $f$ with respect to an additional variable $x_0$ and produce a more geometric counterexample:

\begin{corollary}\label{cor:proj}
Given any $\epsilon>0$ and $n>2$, there exist pairs $(\mathbb{P}_K^n,\Delta)$ such that for an open set $U\in (\mathbb{P}_K^n)^\vee$ and all $H\in U$,
	$$\tau(\mathbb{P}_K^n,\Delta) \neq \tau(\mathbb{P}_K^n,\Delta + \epsilon H)$$
\end{corollary}

Utilizing \autoref{thm:Cex}, one can easily produce explicit counterexamples.

\begin{corollary}[Dimension 4]\label{thm:Ex1}
Let $K=K^p$ be an infinite perfect field of characteristic $p>0$, and let $S = K[x,y,z,w]$. Then there exists $f\in S$ for which a general hyperplane $H=V(l)$ has the property
$$\tau(\mathbb{A}^4,f^\frac{1}{p}) \neq \tau(\mathbb{A}^4,(l\cdot f)^\frac{1}{p}) $$
\end{corollary}
\bpf
	Let $F_*f = F_*(xyz)^{p-1}\left[xF_*w^{p-1} + y^{p-1}zF_*w^{p-2}+y^{p-2}z^2F_*w^{p-3} + \ldots + yz^{p-1}\right]$.  
	
	It immediately meets the conditions of the theorem:
	$$x \notin \langle y^{p-1}z,y^{p-2}z^2,\ldots,yz^{p-1},x^2,xy,xz\rangle_S$$
	$$y^{p-1}z \notin \langle y^{p-2}z^2,\ldots,yz^{p-1}\rangle + \langle y^pz,x^2,xy,xz\rangle_S$$
	$$\vdots$$
	$$yz^{p-1} \notin \langle x,y,z\rangle\cdot \langle x,y^{p-1}z,y^{p-2}z^2,\ldots,yz^{p-1}\rangle_S.$$
	
	\noindent Thus for any $l=c_0+c_xx+c_yy+c_zz+c_ww$ with no $c_\cdot = 0$, one concludes directly that $\tau(\mathbb{A}^4,f^{\frac{1}{p}})\neq \tau(\mathbb{A}^4,(f\cdot l)^{\frac{1}{p}})$. Just for comparison's sake, in this case we are comparing $\langle x, y^{p-1}z,y^{p-2}z^2,\ldots,yz^{p-1}\rangle_S$ with
	$$\langle c_0x+c_wy^{p-1}z, c_0y^{p-1}z + c_wy^{p-2}z^2,\ldots, c_0y^2z^{p-2} + c_wyz^{p-1},c_0yz^{p-1}+c_wwx, $$
	$$ x^2,xy^{p-1}z, xy^{p-2}z^2,\ldots,xyz^{p-1},y^pz, y^{p-1}z^2,\ldots,y^2z^{p-1},yz^p \rangle_S$$
	
	From the monomials it is clear that these generators can not yield a $xw^j$ nor $y^iz^{p-i}w^j$ for any $j>0$ or $i=1,\ldots,p-1$.  Thus, the hyperplanes that satisfy restriction are a subset of the closed subset of hyperplanes with at least one coefficient $0$.  To be more precise, $c_w$ must be $0$ and $c_0\neq 0$ for the equations to work out in such a way that the two ideals agree.\epf
	
\begin{corollary}[Dimension $n\geq 3$]\label{thm:Ex2}
	Consider $S= K[x_1,\ldots,x_n]$ with $K$ infinite perfect of characteristic $p>0$, and $n\geq 3$.  Let $H_0,\ldots,H_{p-1}\subseteq\mathbb{A}^{n-1}$ be general hyperplanes through the origin, viewed as $V(x_n)\subset \mathbb{A}^n$ with $H_i=V(l_i)$ (thus $l_i = c_{i,1}x_1+ \ldots c_{i,n-1}x_{n-1}$ for some $c_{i,j}\in K$). Consider $f_i$ the product of all but the $i^{\text{th}}$ of these hyperplanes:
	$$f_i = \prod_{j=0,\ldots,\hat{i},\ldots, p-1} l_j$$

Then these $f_i$ satisfy the conditions of the theorem.  Thus $F_*f = F_*(x_1\cdots x_{n-1})^{p-1}\sum f_iF_*x_n^i$ yields a $n$-dimensional counterexample in any positive characteristic to Bertini for test ideals.
\end{corollary}	
\bpf
As each $f_i$ is homogeneous of degree $p-1$, the condition $(\star)$ is equivalent to
	$$f_i\notin \langle f_0,\ldots, f_{i-1}\rangle$$
	
\noindent This is arranged by construction, since 
$$V(\langle f_0,\ldots, f_{i-1}\rangle) = \bigcap_{j=0,\ldots,i-1}V(f_j) = \bigcap_{j=0,\ldots,i-1}\left(\bigcup_{k\neq j} H_k\right) = H_i\cup\ldots\cup H_{p-1}$$

\noindent In particular, it contains $H_i$ which $V(f_i)$ does not, so $V(f_i)\not\supseteq V(\langle f_0,\ldots, f_{i-1}\rangle)$.  Therefore, ($\star$) is met (all ideals involved are radical), and one can conclude that for a general hyperplane $H=V(l)$, 
	$$\tau(\mathbb{A}^n,f^\frac{1}{p}) \neq \tau(\mathbb{A}^n,(l\cdot f)^\frac{1}{p}).$$\epf
	
In addition, I show that the direct analog of Bertini's Second Theorem for test ideals fails in general.  Namely, there exists $(X,\Delta)$ with the property that for a general hyperplane $H$, one has
$$\tau(X,\Delta)\cdot\mathscr{O}_H \neq \tau(H,\Delta|_H).$$

\begin{proposition}\label{prop:FailSB1}
	(Dimension 4) Consider \autoref{thm:Ex1}.  Let $D$ be the divisor associated to $f$.  Then for a general choice of $H=V(l)$,
	$$\tau(\mathbb{P}^4,\frac{1}{p}D)|_H \neq \tau(H,\frac{1}{p}D|_H)$$ 
\end{proposition}
\bpf
	Let $l=c_0+c_xx+c_yy+c_zz+c_ww$ with each $c_i\neq 0$.	By \autoref{prop:Comp}, It is enough to show that $x\notin\tau(H,\bar{f}^{\frac{1}{p}})$.  Begin by eliminating $w$ from all equations by virtue of the relation $w = -c_w^{-1}(c_0 + c_xx + c_yy + c_zz)$.  Then 
	$$F_*\bar{f} = \left[(c_w^{-1}c_0)^{\frac{p-1}{p}}x + (c_w^{-1}c_0)^{\frac{p-2}{p}}y^{p-1}z + \ldots + yz^{p-1} \right]F_*(xyz)^{p-1} + \sum_{\boldsymbol{\alpha}<\boldsymbol{p-1}}s_{\boldsymbol{\alpha}}F_*x^{\boldsymbol{\alpha}}$$
Now, utilizing the fact that the test ideal can be computed as in \autoref{prop:Comp}, but with respect to any basis of $F_* (S/\langle l) \rangle$ over $S/\langle l \rangle$, such as the $(xyz)^{\boldsymbol{\alpha}}$ with $\boldsymbol{0}\leq \boldsymbol{\alpha} < \boldsymbol{p}$ basis, it is easily seen that
	$$\tau(H,\bar{f}^{\frac{1}{p}}) = \langle (c_w^{-1}c_0)^{\frac{p-1}{p}}x + (c_w^{-1}c_0)^{\frac{p-2}{p}}y^{p-1}z + \ldots + yz^{p-1},s_{\boldsymbol{\alpha}} \rangle$$	
	It is also easy to see that each $s_{\boldsymbol{\alpha}} \in \langle x,y,z\rangle \cdot \left(\tau(S,f^{\frac{1}{p}}) + \langle l \rangle\right)$, since the replacement of $w$ with $-c_w^{-1}(c_0 + c_xx + c_yy + c_zz)$ yields either a constant coefficient (accounted for in the explicitly written down generator) or some higher multiple of $x$, $y$, or $z$, and further every term of the original $f$ had $F_*(xyz)^{p-1}$ in it already.  Thus, it is clear that $0\neq x\notin \tau(H,\bar{f}^{\frac{1}{p}})$, whereas $x\in \tau(\P^4,f^{\frac{1}{p}})$.\epf
	
Similar considerations apply to \autoref{thm:Ex2}.  However, in the case where the $f_i$ are chosen homogeneous of fixed degree, there is an easy method to detect failure of Bertini:

\begin{theorem}\label{thm:BertiniType}
	Let $S = K[x_1,\ldots,x_n]$, for $K$ a perfect infinite field and let $$F_*f = \sum_{i=0}^{p-1}f_iF_* x_1^{p-1}\cdots x_{n-1}^{p-1} x_n^i$$ with each $f_i\in K[x_1,\ldots,x_{n-1}]$, homogeneous of the same degree $d$, and such that the $f_i$ span a $K$-vectorspace of $S_d$ of dimension at least 2.  Then for a general choice of hyperplane $H$,
$$\tau(S,f^\frac{1}{p})\cdot S/\langle l \rangle \neq \tau(S/\langle l \rangle,\bar{f}^\frac{1}{p})$$
\end{theorem}

\bpf
	Let
$$l=c_0^p+c_1^px_1+\ldots+c_n^px_n.$$
$\tau(S,f^{\frac{1}{p}}) = \langle f_i \ | \ i=0,1,\ldots,p-1\rangle$, as usual.  Assume that $c_n\neq 0$ in the presentation of $l$.  Then $S/\langle l \rangle \cong K[x_1,\ldots,x_{n-1}]$ via
$$x_i\mapsto x_i \ \ \ \ i=1,\ldots,n-1$$
$$x_n \mapsto -c_n^{-p}(c_0^p+c_1^px_1\ldots+c_{n-1}^px_{n-1}).$$
Note that $\bar{f_i}$ are also homogeneous of degree $d$ using the standard grading of $K[x_1,\ldots,x_{n-1}]$. Moreover, one can compute $\tau(S/\langle l \rangle,\bar{f}^\frac{1}{p})$ with respect to any basis of $F_*R$, such as $F_*x_1^{\alpha_1}\cdots x_{n-1}^{\alpha_{n-1}}$, where $0\leq \alpha_j\leq p-1$. So \begin{align*}
F_*\overline{f} = & \sum_{i=0}^{p-1} f_iF_* x_1^{p-1}\cdots x_{n-1}^{p-1} x_n^i\\
= & \sum_{i=0}^{p-1} f_iF_* x_1^{p-1}\cdots x_{n-1}^{p-1}\cdot(-c_n)^{-ip}(c_0^p+c_1^px_1+\ldots+c_{n-1}^px_{n-1})^i.\end{align*}
Therefore, $\bar{f}$ expressed in this basis has non-zero coefficients associated to $F_*\boldsymbol{x^{p-1+\theta}}$, where exponents are taken $\mod p$ and $0\leq |\boldsymbol{\theta}|\leq p-1$.  Applying trace of $S/\langle l\rangle$, the generators come in two types:
\begin{itemize}
	\item $\sum_{i=0}^{p-1} (-c_0c_n^{-1})^i\cdot f_i$, the coefficient of $F_*\boldsymbol{x^{p-1}} = F_*x_1^{p-1}\cdots x_{n-1}^{p-1}$.
	\item $\boldsymbol{x^{\lceil \frac{\theta}{p}\rceil}}\sum_{i=|\boldsymbol{\theta}|}^{p-1} d_{\boldsymbol{\theta},i}\cdot f_i$ the coefficient of $F_*\boldsymbol{x^{p-1+\theta}}$ ($\mod p$ exponents)
\end{itemize}
where $d_{\boldsymbol{\theta},i}=c_0^{i-|\boldsymbol{\theta}|}(-c_n)^{-i}\boldsymbol{c^\theta}\in K$.  It is clear that $\tau(S/\langle l \rangle,\bar{f}^{\frac{1}{p}})$ is homogeneous and
$$\tau(S/\langle l \rangle,\bar{f}^{\frac{1}{p}}) \subseteq \langle \sum_i (c_0c_n^{-1})^if_i\rangle + \mathfrak{p}\cdot \tau(S, f^{\frac{1}{p}})\cdot S/\langle l \rangle$$
where $\mathfrak{p} = \langle x_1,\ldots,x_{n-1} \rangle$.  Intersecting $\tau(S/\langle l \rangle, \bar{f}^\frac{1}{p})$ with $(S/\langle l\rangle)_d$ yields a 1-dimensional vector space $\langle \sum_i (c_0c_n^{-1})^if_i  \rangle_K$.  Therefore, since both ideals are homogeneous, one can compare them degree-wise, and note that 
$$\left(\tau(S,f^\frac{1}{p})\cdot S/\langle l \rangle\right)_d \neq \tau(S/\langle l \rangle,\bar{f}^{\frac{1}{p}})_d.$$
and therefore, the two ideals are distinct.  This completes the proof.
\epf	

\begin{corollary}\label{prop:FailSB2}
(Dimension $n\geq 3$) Consider the situation of \autoref{thm:Ex2}.  If $D$ is the divisor associated to $f$, then for a general choice of $H=V(l)$,
$$\tau(\mathbb{P}^n,\frac{1}{p}D)|_H \neq \tau(H,\frac{1}{p}D|_H)$$ 
\end{corollary}

\bpf
Choosing each of the $l_i$ as in \autoref{thm:Ex2} implies directly that the $f_i$ are homogeneous of degree $p-1$.  Therefore, since the $l_i$ are chosen generally, any two of the $f_i$ are linearly independent, and thus the condition of \autoref{thm:BertiniType} is satisfied.
\epf

This answers \cite[Question 8.3]{SZ} in the negative in any dimension larger than 2.  Moreover, a similar argument provides a counterexample to a natural extension of a question of Hochster and Huneke.

\begin{theorem}
\cite[Theorem 7.3(c)]{HH} Let $\varphi:R\rightarrow S$ be a flat homomorphism of Noetherian rings in characteristic $p>0$.  Suppose that $R$ is F-regular, S is excellent, and that $\varphi$ has regular fibers.  Then $S$ is also F-regular.
\end{theorem}

This brings about the following question: 

\begin{Question}If $\pi: X\rightarrow Y$ is a flat morphism of characteristic $p>0$ \ Noetherian schemes with regular fibers, with $\Delta$ a $\Q$-Cartier divisor on $Y$, is it true that $\tau(Y,\Delta_Y)\cdot\OO_X = \tau(X,\pi^*\Delta_Y)$?
\end{Question}

  The result of Hochster and Huneke holds for Strongly F-regular pairs, where both test ideals are their respective sheaves of rings.  However, \autoref{prop:FailSB1} (or similarly \autoref{prop:FailSB2}) can be used to show that this question is false more generally.

\begin{theorem}\label{thm:HHfalse}
Let $K$ be a perfect field, and let $Z$ be the reduced closed subscheme of $\P^n\times_K(\P^n)^\vee$ defined to be the closure of the set $\{(P,H) \ | \ P\in H\}$.  Then $Z\rightarrow \P^n$ is a flat morphism, with regular fibers, such that there exists $\Delta$ a $\Q$-divisor on $\P^n$ with the property that 
$$\tau(\P^n,\Delta) \cdot \O_Z \neq \tau(Z,\pi^*\Delta).$$
\end{theorem}

\bpf
The statement that  $Z\rightarrow \PP^n$ is a flat morphism is proved in the course of \cite[Theorem 3.7]{SZ}.  Let $\Delta = \frac{1}{p}\text{div}(f)$, where $f$ is chosen to satisfy the assumptions of \autoref{thm:BertiniType}.  Note that to show that the two ideals of $\mathcal{O}_Z$ differ, it suffices to show that their localizations differ.  Let $P\in \P^n$, and consider the maximal ideal corresponding to a point $\mathfrak{m} = (P,H)\in Z$.  Then the localization map has the form
$$K[x_1,\ldots,x_n]_{\langle \underline{x}\rangle} = \mathcal{O}_{\P^n,P} \rightarrow \mathcal{O}_{Z,(P,H)} \cong K[x_1,\ldots,x_n,y_1,\ldots,y_n]/\langle y_1x_1+\ldots+y_nx_n\rangle_{\langle \underline{x},\underline{y}\rangle}.$$
Here the underline is representing $\underline{x} = x_1,\ldots,x_n$, and similarly for $\underline{y}$.  Note this is the case since every hyperplane passing through the origin necessarily has this form. Since test ideals localize, we can perform a further localization at the prime ideal ${\langle \underline{x}\rangle \subseteq \mathcal{O}_{Z,(P,H)}}$.  It suffices to prove that
$$\tau(\P^n,\Delta)\cdot \mathcal{O}_{Z,\langle \underline{x}\rangle}=\tau(\P^n,\Delta)_{\langle \underline{x}\rangle}\cdot \mathcal{O}_{Z,\langle \underline{x}\rangle} \neq \tau(Z,\pi^*\Delta)_{\langle \underline{x}\rangle}$$

The left hand side is simply $\tau(S,f^{\frac{1}{p}})_{\langle \underline{x}\rangle} = \langle f_i \ | \ i=0,1,\ldots,p-1\rangle\cdot S_{\langle \underline{x}\rangle}$.  The right hand side is $$\tau\left(k(y_1,\ldots,y_n)[x_1,\ldots,x_n]/\langle y_1x_1+\ldots+y_nx_n\rangle_{\langle \underline{x}\rangle},f^{\frac{1}{p}}\right).$$
Call this ring $R$ for simplicity.  This ideal can be computed by applying the trace operator $\Phi_R$ to the ideal $\langle f\rangle_R$.  Similar to the proof of \autoref{thm:BertiniType}, we can compute the test ideal with respect to any basis of $F_*R$ over $R$.  Noting that
$$x_n = y_n^{-1}\left(y_1x_1+\ldots+y_{n-1}x_{n-1} \right)$$
the basis of interest is
$$\left\lbrace F_*\left(y_1^{\beta_1}\cdots y_n^{\beta_n}\cdot x_1^{\alpha_1}\cdots x_{n-1}^{\alpha_{n-1}}\right) \ : \ 0\leq \alpha_i,\beta_i < p\right\rbrace$$
Therefore, ones concludes $\tau(Z,\pi^*\Delta)_{\langle \underline{x}\rangle}$ is homogeneous in $x_1,\ldots,x_{n-1}$ and has at most one generator of $x$-degree $d$, whereas $\tau(\P^n,\Delta)$ has at least $2$ by the assumptions imposed on $f$.  This completes the proof.
\epf
As a final consideration, we ask about the Bertini statement in dimension 2.  It is immediately clear that \autoref{thm:Cex} cannot be directly extended to dimension 2, as finding a collection of $f_i$ meeting $(\star)$ is clearly impossible using only a single variable.  In fact, some simple geometry allows us to conclude that it is in fact true.

\begin{theorem}\label{thm:dim2}
Let $X\subseteq \PP^n$ be a normal projective surface over an infinite perfect field $K$ of characteristic $p>0$, and let $\Delta\geq 0$ be an effective $\Q$-divisor, with $K_X+\Delta$ $\Q$-Cartier.  Then for a general hyperplane $H\in (\P_K^n)^\vee$,
$$ \tau(X,\Delta)\cdot \OO_H = \tau(X\cap H,\Delta|_H)$$

Furthermore, for every $0\leq \epsilon < 1$, 
$$ \tau(X,\Delta) = \tau(X,\Delta + \epsilon H)$$
\end{theorem}

\bpf
	Let $\Sigma$ be the locus containing singularities of $X$ and singularities of $\Delta$.  Then by normality, $\Sigma$ is a finite set of closed points.  Take $H$ any hyperplane not intersecting $\Sigma$, with $X\cap H$ regular.  Then by completing the local ring at a point of $x\in \Delta\cap H$, by regularity Cohen structure theorem implies $\OO_{X,x} = K\llbracket x,y\rrbracket$, and by a linear change of coordinates, we can assume $H=V(y)$ and $\Delta = t\cdot V(x)$ for some positive rational number $t$.  In this complete local setting, $\OO_H = K\llbracket x \rrbracket$, and $\Delta|_H$ is simply $t\cdot V(x)$.  Finally, the following computations follow directly from \autoref{Sect:Comp}:	
	$$\tau(K\llbracket x,y\rrbracket,x^t) \otimes_{K\llbracket x,y\rrbracket} K\llbracket x\rrbracket = \langle x^{\lfloor t \rfloor}\rangle = \tau(K\llbracket x\rrbracket,x^t)$$
	
	Finally, since every point of $X\cap H$ and $\Delta \cap H$ is smooth in $X,\Delta$ respectively, the test ideal sheaves agree as well.
	
	The second statement follows by the same logic, as $\tau(X,\Delta + \epsilon H) = \langle x^{\lfloor t \rfloor} \rangle$ in $K\llbracket x,y\rrbracket$.\epf
	
\section{Further Questions}\label{Sect:Questions}

The previous section implies that an analog of Bertini's Theorem for test ideals is not possible for pairs $(X,\Delta)$ in general.  However, there are also classes for which such a theorem is possible.  For example, if $f$ is of the form $f_\alpha F_*x^\alpha$, or if the $F_*x^\alpha$ are in sufficiently low degrees, or even in many cases where the $f_i$ are not quite as independent as condition $(\star)$ enforces.  A more general class comes from \autoref{thm:SZ}, which states that for Sharply F-pure or Strongly F-regular pairs $(X,\Delta)$, an analog of the second theorem of Bertini in fact holds.  So a general question one could pose would be which pairs $(X,\Delta)$ have such a property?  Or more specifically, is there some geometric or arithmetic property that is governing whether or not the theorem holds?

  One particular thing to note is that \autoref{thm:Ex2} has $\Delta$ of a fairly large degree; $(p+n)\cdot(p-1)$ where $n$ is the dimension.  It is possible that the theorem holds up to a particular degree for a given dimension of variety.  One could also ask if there is some type of locus determined by the test ideal for which a hyperplane not intersecting the said locus is sufficient to guarantee that the restriction statement holds, and what properties this locus has.  In particular, if it was a finite collection of points, a general hyperplane would always miss this loci and the Bertini theorem would hold.

Two geometric properties of interest in positive characteristic Algebraic Geometry and Commutative Algebra, which are related to the test ideal, are F-injectivity and F-rationality.  The direct analog of the second theorem of Bertini for F-injectivity is known to be false \cite[Proposition 7.4]{SZ} (still open for normal F-injective, as the counter-example used the Weakly normal construction of \cite{CGM}), however it is currently unknown whether the same is true for F-rationality.

Finally, a very important question that remains open is whether or not one can use this type of technique to find counterexamples to the original Hochster-Huneke question.  Namely, could it still be true when considering a variety instead of pairs?

\begin{Question}
	For an irreducible projective variety $X$, and a very ample linear series $|D|$, is it true that for a general choice of $\Delta\in |D|$, one has that $\tau(X) = \tau(\Delta\cap X)$?
\end{Question}

The answer is very likely no based on the findings of this paper, but an explicit counterexample is still unknown to the author.  In particular, the class of $f$ provided from \autoref{thm:Cex}, $X=V(f)$, for which the trace map has a particularly nice form $\Phi_X(g) = \Phi(g F_*f^{p-1})$, cannot be used because none of these $X$ are irreducible varieties.

\bibliographystyle{skalpha}
\bibliography{BibFile}
\end{document}